\documentclass[12pt,leqno]{article}

\usepackage{latexsym,amsmath,amsfonts}

\newtheorem{theorem}{Theorem}[section]
\newtheorem{lemma}[theorem]{Lemma}
\newtheorem{proposition}[theorem]{Proposition}

\newtheorem{remark}[theorem]{Remark}
\newtheorem{example}[theorem]{Example}
\newtheorem{definition}[theorem]{Definition}

\def\RN{\mathbb R^N}
\def\rh{\rightharpoonup}

\def\pf{\noindent\textbf{Proof}\ \ }

\title{\bf Ground state for the Schr\"{o}dinger operator with \\  the weighted Hardy potential
\author
{J. Chabrowski\\
{\small Department of Mathematics}\\
{\small The University of Queensland}\\
{\small St. Lucia 4072, Qld, Australia}\\
e-mail:jhc@maths.uq.edu.au\\ \\
K. Tintarev\\
Department of Mathematics\\
Uppsala University P.O. box 480\\
751 06 Uppsala, Sweden\\
e-mail: tintarev@math.uu.se}}

\date{}

\begin{document}

\baselineskip15pt \maketitle
\renewcommand{\thefootnote}{}

 \begin{abstract}
 We establish the existence of ground states on $\mathbb R^N$ for the
 Laplace operator involving the  Hardy type potential. This gives rise to the existence of the principal eigenfunctions
 for the Laplace operator involving weighted Hardy potentials.
  We also obtain a higher integrability property for the principal eigenfunction. This is used to examine  the behaviour
  of the principal eigenfunction around $0$.
 \end{abstract}

\section{Introduction}
\renewcommand{\theequation}{1.\arabic{equation}}
\setcounter{equation}{0}

\medskip

In this paper we investigate the existence of ground states of the  Schr\"{o}dinger operator associated with
the quadratic form
\begin{equation}\label{a1}
Q_V(u)=\int_{\RN} \bigl(|\nabla u|^2-\Lambda_VV(x)u^2\bigr)\,dx, \,\,  u \in C_{\circ}^\infty(\RN), \,\, N\geq
3,
\end{equation}
where $V$ belongs to the Lorentz space $L^{\frac N2,\infty}(\RN)$ and  $\Lambda_V$ is the largest constant
(whenever exists) for which the form $Q_V$ is nonnegative. This assumption implies that the potential term
$\int_{\RN} V(x)u^2\,dx$ is continuous in $D^{1,2}(\RN)$, where $D^{1,2}(\RN)$ is the Sobolev space obtained
as the completion of $C_{\circ}^\infty(\RN)$ with respect to the norm
\[
\|u\|_{D^{1,2}}^2=\int_{\RN} |\nabla u|^2\,dx.
\]
We are mainly interested in the case of the Hardy type potential $V(x)=\frac {m(x)}{|x|^2}$ with $m \in
L^{\infty}(\RN)$. Assuming that $V$ is positive on a set of positive measure, the constant $\Lambda_V$ is
given by the variational problem
\begin{equation}\label{a2}
\Lambda_V=\inf_{u \in D^{1,2}(\RN), \, \int_{\RN} Vu^2\,dx=1}\int_{\RN} |\nabla u|^2\,dx
 \end{equation}
 and  the continuity of $\int_{\RN}V(x)u^2\,dx$  implies that $\Lambda_V>0$. If problem (\ref{a2}) has a
 minimizer $u$, then it satisfies the equation
 \begin{equation}\label{a3}
 -\Delta u-\Lambda_V V(x)u=0.
 \end{equation}
 A solution of (\ref{a3})  is understood in the weak sense
\begin{equation}\label{a4}
\int_{\RN} \nabla u\nabla \phi \, dx=\Lambda_V \int_{\RN} V(x) u\phi\,dx,
\end{equation}
for every $\phi \in D^{1,2}(\RN)$.

Since $|u|$ is also a minimizer for $\Lambda_V$, we may assume that $u \geq 0$ a.e. on $\RN$. In particular,
when $V(x)=\frac {m(x)}{|x|^2}$ with $m \in L^{\infty}(\RN)$, then $u>0$ on $\RN$ by the Harnack inequality
\cite{GT}. If the potential term is weakly continuous in $D^{1,2}(\RN)$, for example, when $V(x)=\frac
{m(x)}{|x|^2}$ with $m \in L^{\infty}(\RN)$ and $\lim_{|x| \to \infty} m(x)=\lim_{x \to 0} m(x)=0$, then there
exists a minimizer for $\Lambda_V$. We will call the minimizer of (\ref{a2}) a {\em ground state of finite energy}.
In general, (\ref{a2}) may not have a minimizer. This is the case for the Hardy potential $V(x)=\frac
1{|x|^2}$ with the corresponding optimal constant $\Lambda_V=\Lambda_N=\biggl(\frac {N-2}2\biggr)^2$.  In
fact, the ground state of finite energy is a particular case of the generalized ground state, defined as
follows (see \cite{MU}, \cite{PT1} and \cite{PT2}).

\medskip
\begin{definition}\label{d1}
Let $\Omega \subset \RN$ be an open set, and let $Q_V$ be as in (\ref{a1}).  A sequence of nonnegative
functions $v_k \in C_{\circ}^\infty(\Omega)$ is said to be a {\it null-sequence }  for the functional $Q_V$ if
$Q_V(v_k) \to 0$, as $k \to \infty$,  and there exists a nonnegative function $\psi \in
C_{\circ}^\infty(\Omega)$ such that $\int_{\Omega}\psi v_k\,dx=1$ for each $k$.
\end{definition}

\medskip

Let us recall that the capacity   of a compact set $E$ relative to an open set $\Omega \subset \RN$, with
$E \subset \Omega$, is given by
\[
{\mbox{\rm cap }}(E,\Omega)=\inf\{ \int_{\Omega} |\nabla u|^2\,dx;\, u \in C_{\circ}^{\infty}(\Omega), \,
\mbox{\rm with }\, u(x) \geq 1 \,\mbox{ on }\, E\}.
\]
In the case $\Omega=\RN$ we use notation $\mbox{\rm cap }(E)$ (see \cite{MA}).

\medskip

We can  now formulate  the following "ground state alternative" (see \cite{PT1}, \cite{PT2}).

\medskip

\begin{theorem}\label{t1}
Let $V$ be a measurable function bounded on every compact subset of $\Omega =\RN-Z$, where $Z$ is a closed set of
capacity zero, and assume that $Q_V(u)\geq 0$ for all $u \in C_{\circ}^\infty(\Omega)$. Then, if $Q_V$ admits
a null sequence $v_k$, then the sequence $v_k$ converges weakly in $H^1_{\mbox{\rm loc}}(\RN)$ to a unique (up
to a multiplicative constant ) positive solution of (\ref{a3}).
\end{theorem}

\medskip

This theorem gives rise to the definition of the generalized ground state.

\begin{definition}\label{d2}
A unique positive solution $v$ of (\ref{a3}) is called a generalized ground state of the functional $Q_V$, if
the functional admits a null sequence weakly convergent to $v$.
\end{definition}

\medskip

If $V(x)=\frac 1{|x|^2}$, the functional $Q_V$ has a ground state $v(x)=|x|^{\frac {2-N}2}$ of infinite
$D^{1,2}$ norm, while (\ref{a2}) has no minimizer in $D^{1,2}(\RN)$.

\medskip

It is important to note that the functional $Q_V$ with the optimal constant $\Lambda_V$ does not necessarily
have a ground state. We quote the following statement from \cite{PT2}.

\medskip

\begin{theorem}\label{t2}
Let $V$ be a measurable function bounded on every compact subset of $\Omega =\RN-Z$, where $Z$ is a closed set
of capacity zero, and assume that $Q_V(u)\geq 0$ for all $u \in C_{\circ}^\infty(\Omega)$. Then either $Q_V$
admits a null sequence, or there exists a function $W$, positive and continuous on $\Omega$, such that
\begin{equation}\label{a5}
Q_V(u) \geq \int_{\RN} W(x)u^2\,dx.
\end{equation}
\end{theorem}

\medskip

For example, let $m$ be a continuous function on $\RN-\{0\}$ such that $m(x)=\frac 1{|x|^2}$ for $0<|x| \leq
1$, $m(x) \in [\frac 12,1]$ for $|x| \in (1,2)$ and $m(x)=\frac{1}{2|x|^2}$ for $|x| \geq 2$. Then
$\Lambda_V=\biggl(\frac {N-2}2\biggr)^2$ and the functional $Q_V$ does not admit a null sequence.  From
Theorem~\ref{t2} follows that $Q_V$ satisfies (\ref{a5}) with some function $W$ positive on $\RN-\{0\}$.

\medskip

Obviously, ground states of finite $D^{1,2}$ norm are principal eigenfunctions of (\ref{a3}). There is a quite
extensive literature on principal eigenfunctions with indefinite weight functions for elliptic operators on
$\RN$, or on unbounded domains of $\RN$, with the Dirichlet boundary conditions. We mention papers \cite{AL},
\cite{BCF}, \cite{BT}, \cite{JI}, \cite{LY}, \cite{MU}, \cite{RH}, \cite{SM}, \cite{SW}, where the existence
of principal eigenfunctions has been established under various assumptions on weight functions. These
conditions require that a potential belongs to some Lebesgue space, for example $L^p(\RN)$ with $p >\frac N2$.
These results have been recently greatly improved in papers \cite{ALR} and \cite{VI}, where potentials from
the Lorentz spaces have been considered. To describe the results from \cite{ALR} and \cite{VI} we recall the
definition of the Lorentz space \cite{BL}, \cite{HU}, \cite{LO}.

\medskip

Let $f:\, \RN \to \mathbb R$ be a measurable function. We define the distribution function $\alpha_f$ and a
nonincreasing  rearrangement $f^*$ of $f$ in the following way
\[
\alpha_f(s)=|\{x \in \RN;\, |f(x)|>s\} \,\, \mbox{ and }\,\, f^*(t)=\inf\{s>0;\, \alpha_f(s) \leq t\}.
\]
We now set
\[
\|f\|_{(p,q)}^*=\left\{
\begin{array}{ll}
\biggl(\int_0^\infty [t^{\frac 1p}f^*(t)]^q\frac {dt}t\biggr)^{\frac 1q}& \,\,\, \mbox{ if } \, 1\leq p,q <\infty,\\
\sup_{t>0} t^{\frac 1p} f^*(t) &\,\,\, \mbox{ if }\, 1\leq p \leq \infty, q=\infty.
\end{array}
\right.
\]
The Lorentz space $L^{p,q}(\RN)$ is defined by
\[
L^{p,q}(\RN)=\{f \in L_{\mbox{ loc }}^1(\RN);\,\|f\|_{(p,q)}^*<\infty\}.
\]
The functional $\|f\|_{(p,q)}^*$ is only a quasi-norm. To obtain a norm we replace $f$ by $f^{**}(t)=\frac 1t
\int_0^tf^*(s)\,dx$ in the definition of $\|f\|_{(p,q)}^*$, that is, the norm is given by
\[
\|f\|_{(p,q)}=\left\{
\begin{array}{ll}
\biggl(\int_0^\infty [t^{\frac 1p}f^{**}(t)]^q\frac {dt}t\biggr)^{\frac 1q}& \,\,\, \mbox{ if } \, 1\leq p,q <\infty,\\
\sup_{t>0} t^{\frac 1p} f^{**}(t) &\,\,\, \mbox{ if }\, 1\leq p \leq \infty, q=\infty.
\end{array}
\right.
\]
$L^{p,q}(\RN)$ equipped with the norm $\|f\|_{(p,q)}$ is a Banach space.

\medskip

In paper \cite{VI} the existence of principal eigenfunctions has been established for weights belonging to
$\bigcup_{1 \leq q <\infty} L^{\frac N2,q}(\RN)$. This was extended in \cite{ALR} to a larger class of weights
${\mathcal{F}}_{\frac N2}$ obtained as the completion of $C_{\circ}^{\infty}(\RN)$ in norm $\|\cdot\|_{\frac
N2,\infty}$.

\medskip

However, these conditions do not cover the singular weight functions considered in this paper. By contrast, in
our approach we give an exact upper bound for the principal eigenvalue which allows us to prove the existence
of the principal eigenfunction. We point out that if $V \in L^{\frac N2,\infty}(\RN)$, then the functional
$\int_{\RN} V(x)u^2\,dx$ is continuous on $D^{1,2}(\RN)$, but not necessarily weakly continuous.

\medskip

The paper is organized as follows. In Section 2 we prove the existence of  minimizers with finite norm
$D^{1,2}(\RN)$ and also with infinite norm $D^{1,2}(\RN)$. In Section 3  we discuss perturbation of a given
quadratic form $Q_{V_{\circ}}$ with $V_{\circ}\in L^{\frac N2,\infty}(\RN)$. We show that if $Q_{V_{\circ}}$
has ground state, then this property is stable under small perturbations of $V_{\circ}$. This is not true if
$Q_{V_{\circ}}$ does not have a ground state; rather it is  stable under larger perturbation of $V_{\circ}$.
The final Section is devoted to a higher integrability property of  minimizers of $Q_{V_{\circ}}$ in the case
where $V_{\circ}(x)=\frac {m(x)}{|x|^2}$ with $m \in L^{\infty}(\RN)$. We also examine the behaviour of the
principal eigenfunction around $0$.

\medskip

Throughout this paper, in a given Banach space we denote strong convergence by $"\to"$ and weak convergence by
$"\rh"$. The norms in the Lebesgue space $L^p(\Omega)$, $1\leq p\leq \infty$, are denoted by $\|u\|_p$.

\medskip

\section{Existence of minimizers}
\renewcommand{\theequation}{2.\arabic{equation}}
\setcounter{equation}{0}

\medskip

We consider the Hardy type potential $V(x)=\frac {m(x)}{|x|^2}$ with $m \in L^{\infty}(\RN)$. In
Theorem~\ref{t3} we formulate conditions on $m$ guaranteeing the existence of a principal eigenfunction. Let
$\gamma_+>1$ and $\gamma_->1$. In our approach to problem (\ref{a2}) the following two limits play an
important role: it is assumed  that the following limits exist a.e.
\begin{equation}\label{a6}
m_+(x)=\lim_{j \in \mathbb N, j \to \infty} m(\gamma_+^jx)
\end{equation}
and
\begin{equation}\label{a7}
m_-(x)=\lim_{j \in \mathbb N, j \to \infty} m(\gamma_-^{-j}x).
\end{equation}
Both functions $m_{\pm}$ satisfy $m_{\pm}(\gamma_{\pm}x)=m_{\pm}(x)$, that is, $m_{\pm}$ are  homogeneous of
degree $0$. We now define the following infima:
\begin{equation}\label{a8}
\Lambda_m=\inf_{u \in D^{1,2}(\RN)-\{0\}} \frac {\int_{\RN} |\nabla u|^2\,dx}{\int_{\RN} \frac
{m(x)}{|x|^2}u^2\,dx},
\end{equation}
(we use the notation $\Lambda_m$ instead of $\Lambda_V$) and
\begin{equation}\label{a9}
\Lambda_{\pm}=\inf_{u \in D^{1,2}(\RN)-\{0\}} \frac {\int_{\RN} |\nabla u|^2\,dx}{\int_{\RN} \frac
{m_{\pm}(x)}{|x|^2}u^2\,dx}.
\end{equation}

\medskip

\begin{lemma}\label{l1}
The following holds true
\begin{equation}\label{a10}
\Lambda_m \leq \min(\Lambda_+,\Lambda_-).
\end{equation}
\end{lemma}
\pf Let $u \in D^{1,2}(\RN)-\{0\}$. Testing $\Lambda_m$ with $\gamma_+^{-\frac {N-2}2}u(\gamma_+^{-j}x)$ gives
\[
\Lambda_m \leq  \frac {\int_{\RN} |\nabla u|^2\,dx}{\int_{\RN} \frac {m(\gamma_+^jx)}{|x|^2}u^2\,dx}.
\]
Letting $j \to \infty$ and using the Lebesgue dominated convergence theorem, we obtain
\[
\Lambda_m \leq \frac {\int_{\RN} |\nabla u|^2\,dx}{\int_{\RN} \frac {m_+(x)}{|x|^2}u^2\,dx}.
\]
The inequality $\Lambda_m \leq \Lambda_+$ follows. The proof of the inequality $\Lambda_m \leq \Lambda_-$ is
similar. $\hfill\Box$

\medskip

In the case when the inequality (\ref{a9}) is strict problem (\ref{a7}) has a minimizer.

\medskip

\begin{theorem}\label{t3}
Assume that the convergence in (\ref{a6}) is uniform on sets $\{x \in \RN; \, |x|\geq R\}$ for every $R>0$ and
that the convergence in (\ref{a7}) is uniform on sets $\{x \in \RN;\, |x| \leq \rho\}$ for every $\rho>0$. If
$\Lambda_m <\min(\Lambda_+,\Lambda_+)$, then problem (\ref{a8})  has a minimizer.
\end{theorem}
\pf Let $\{u_k \} \subset D^{1,2}(\RN)$ be a minimizing sequence for $\Lambda_m$, that is,
\[
\int_{\RN} |\nabla u_k|^2\,dx\to \Lambda_m\,\, \mbox{ and }\,\, \int_{\RN} \frac {m(x)}{|x|^2}u_k^2\,dx=1.
\]
We can assume, up to a subsequence, that $u_k \rh w$ in $D^{1,2}(\RN)$, $L^2(\RN,\frac {dx}{|x|^2})$  and $u_k
\to w$ in $L_{\mbox{ \rm loc}}^2(\RN)$ for some $w \in D^{1,2}(\RN)$. Let $v_k=u_k-w$. We then have
\begin{equation}\label{a11}
1=\int_{\RN} \frac {m(x)}{|x|^2}u_k^2\,dx=\int_{\RN} \frac {m(x)}{|x|^2}w^2\,dx+\int_{\RN} \frac
{m(x)}{|x|^2}v_k^2\,dx+o(1)
\end{equation}
and
\begin{equation}\label{a12}
\Lambda_m=\int_{\RN} |\nabla u_k|^2\,dx+o(1)=\int_{\RN} |\nabla w|^2\,dx+\int_{\RN} |\nabla v_k|^2\,dx+o(1).
\end{equation}
We define a radial function $\chi_+^j \in C^1(\RN)$ such that $0 \leq \chi_+^j(x) \leq 1$, $\chi_+^j(x)=0$ for
$|x| \leq \gamma_-^{-2j}$ and $\chi_+^j(x)=1$  for $|x|>\gamma_+^{2j}$. Let $\chi_-^j(x)=1-\chi_+^j(x)$. In
what follows we use $o^{(j)}_{k\to\infty}(1)$ to denote a quantity such that for each $j \in \mathbb N$, 
$o^{(j)}_{k\to\infty}(1) \to 0$ as $k \to \infty$. Thus
\begin{eqnarray}\label{a13}
\int_{\RN} \frac {m(x)}{|x|^2}v_k^2\,dx&=&\int_{\RN} \frac
{m(x)}{|x|^2}\bigl(v_k\chi_-^j\bigr)^2\,dx+\int_{\RN}
\frac {m(x)}{|x|^2}\bigl(v_k\chi_+^j\bigr)^2\,dx+o^{(j)}_{k\to\infty}(1)\\
&=&\int_{\RN} \frac {m(\gamma_-^{-j}x)}{|x|^2}\bigl(v_k^-\bigr)^2\,dx+\int_{\RN} \frac
{m(\gamma_+^jx)}{|x|^2}\bigl(v_k^+\bigr)^2\,dx+o^{(j)}_{k\to\infty}(1), \nonumber
\end{eqnarray}
where
\[
v_k^-(x)=\gamma_-^{-\frac {N-2}2j}v_k\bigl(\gamma_-^{-j}x\bigr)\chi_-\bigl(\gamma_-^{-j}x\bigr)
\]
and
\[
v_k^+(x)=\gamma_+^{\frac {N-2}2j}v_k\bigl(\gamma_+^jx\bigr)\chi_+\bigl(\gamma_+^jx\bigr).
\]
We now estimate the integrals involving $v_k^-$ and $v_k^+$. We have
\begin{eqnarray*}
\left| \int_{\RN} \frac {m\bigl(\gamma_-^{-j}x\bigr)}{|x|^2}(v_k^-)^2\,dx\right.&-&\left.\int_{\RN} \frac
{m_-(x)}{|x|^2}(v_k^-)^2\,dx\right|
 \leq\left| \int_{|x|<\gamma_-^{-j}} \frac
{m\bigl(\gamma_-^{-j}x\bigr)-m_-(x)}{|x|^2}(v_k^-)^2\,dx\right| \\
&+&\left|\int_{\gamma_-^{-j}<|x|<\gamma_-^j\gamma_+^{2j}} \frac
{m(\gamma_-^{-j}x)-m_-(x)}{|x|^2}(v_k^-)^2\,dx\right| =J_1+J_2.
\end{eqnarray*}
By the uniform convergence  of $m\bigl(\gamma_-^{-j}x\bigr)$ to $m_-(x)$ we see that $J_1 \leq \epsilon$ for
$j$ sufficiently large uniformly in $k$. For $J_2$ we have
\[
J_2 \leq 2\|m\|_{\infty} \int_{\gamma_-^{-2j}<|x|<\gamma_+^{2j}} \frac {v_k^2}{|x|^2}\,dx.
\]
It is clear that $J_2$ is a quantity of type $o^{(j)}_{k\to\infty}(1)$. Therefore, we have
\begin{equation}\label{a14}
\left|\int_{\RN}\frac {m\bigl(\gamma_-^{-j}x\bigr)}{|x|^2}(v_k^-)^2\,dx-\int_{\RN} \frac
{m_-(x)}{|x|^2}(v_k^-)^2\,dx\right|\leq \epsilon+o_j(1).
\end{equation}
In a similar way we obtain
\begin{equation}\label{a15}
\left|\int_{\RN}\frac {m\bigl(\gamma_+^jx\bigr)}{|x|^2}(v_k^+)^2\,dx-\int_{\RN} \frac
{m_+(x)}{|x|^2}(v_k^+)^2\,dx\right|\leq \o^{(j)}_{k\to\infty}(1).
\end{equation}
for $j$ sufficiently large. We now fix $j \in \mathbb N$ so that (\ref{a14}) and (\ref{a15}) hold.
Consequently, we have
\begin{equation}\label{a16}
1\leq \int_{\RN} \frac {m(x)}{|x|^2}w^2\,dx+\int_{\RN} \frac {m_-(x)}{|x|^2}(v_k^-)^2\,dx+\int_{\RN} \frac
{m_+(x)}{|x|^2}(v_k^+)^2\,dx+2\epsilon +o^{(j)}_{k\to\infty}(1).
\end{equation}
We now estimate $\int_{\RN} |\nabla v_k|^2\,dx$ in the following way
\begin{eqnarray*}
\int_{\RN} |\nabla v_k|^2\,dx&=&\int_{\RN} |\nabla (v_k\chi_-^j+v_k\chi_+^j)|^2\,dx\\
&=&\int_{\RN} |\nabla (v_k\chi_-^j)|^2\,dx+\int_{\RN} |\nabla (v_k\chi_+^j)|^2\,dx\\
&+&2\int_{\RN} \nabla(v_k\chi_-^j)\nabla(v_k\chi_+^j)\,dx\\
&=&\int_{\RN} |\nabla v_k^-|^2\,dx+\int_{\RN} |\nabla v_k^+|^2\,dx+2\int_{\RN} |\nabla
v_k|^2\chi_-^j\chi_+^j\,dx\\
&+&2\int_{\RN} v_k\nabla v_k\nabla \chi_-^j\chi_+^j\,dx+2\int_{\RN} v_k\nabla v_k\chi_-^j\nabla \chi_+^j\,dx\\
&+&2\int_{\RN} v_k^2\nabla \chi_-^j\nabla \chi_+^j\,dx\\
&\geq&\int_{\RN} |\nabla v_k^-|^2\,dx+\int_{\RN} |\nabla v_k^+|^2\,dx+2\int_{\RN} v_k\nabla v_k\nabla
\chi_-^j\chi_+^j\,dx\\
&+&2\int_{\RN} v_k\nabla v_k\chi_-^j\nabla \chi_+^j\,dx\\
&+&2\int_{\RN} v_k^2\nabla \chi_-^j\nabla \chi_+^j\,dx.
\end{eqnarray*}
Since $v_k \to 0$ in $L_{\mbox{\rm loc}}^2(\RN)$ we obtain the following estimate
\[
\int_{\RN} |\nabla v_k|^2\,dx\geq \int_{\RN} |\nabla v_k^-|^2\,dx+\int_{\RN} |\nabla v_k^+|^2\,dx+o^{(j)}_{k\to\infty}(1).
\]
This, combined with (\ref{a11}), gives the following estimate
\begin{eqnarray}\label{a17}
\Lambda_m &\geq&\int_{\RN} |\nabla w|^2\,dx+\int_{\RN} |\nabla v_k^-|^2\,dx+\int_{\RN} |\nabla
v_k^+|^2\,dx+o_j(1)\\
&\geq&\Lambda_m \int_{\RN} \frac {m(x)}{|x|^2}w^2\,dx+\Lambda_-\int_{\RN} \frac
{m_-(x)}{|x|^2}(v_k^-)^2\,dx\nonumber\\
&+& \Lambda_+\int_{\RN} \frac {m_+(x)}{|x|^2}(v_k^+)^2\,dx+o^{(j)}_{k\to\infty}(1). \nonumber
\end{eqnarray}
Let $\Lambda_*=\min\bigl(\Lambda_-,\Lambda_+\bigr)$. We deduce from (\ref{a16}) and (\ref{a17}) that
\[
\bigl(\Lambda_*-\Lambda_m\bigr)\biggl(\int_{\RN} \frac {m_-(x)}{|x|^2}(v_k^-)^2\,dx+\frac
{m_+(x)}{|x|^2}(v_k^+)^2\,dx\biggr)\leq 2\epsilon \Lambda_m+o^{(j)}_{k\to\infty}(1).
\]
Letting $k \to \infty$ we obtain
\[
\limsup_{k \to \infty}\biggl(\int_{\RN} \frac {m_-(x)}{|x|^2}(v_k^-)^2\,dx+\frac
{m_+(x)}{|x|^2}(v_k^+)^2\,dx\biggr)\leq \frac {2\epsilon\Lambda_m}{\bigl(\Lambda_*-\Lambda_m\bigr)}.
\]
It then follows from (\ref{a16}) that
\[
1 \leq \int_{\RN} \frac {m(x)}{|x|^2}w^2\,dx+\frac {2\epsilon\Lambda_m}{\bigl(\Lambda_*-\Lambda_m\bigr)}.
\]
Since $\epsilon>0$ is arbitrary we get $\int_{\RN} \frac {m(x)}{|x|^2}w^2\,dx=1$ and the result follows.
$\hfill\Box$

\medskip

In what follows, we use denote by $m(\infty)=\lim_{|x| \to \infty} m(x)$, assuming that this limit exists.  As
a direct consequence of Theorem~\ref{t3} we obtain the following result.

\medskip

\begin{theorem}\label{t4}
Let $m \in L^{\infty}(\RN)$ and assume that $m$ is continuous at $0$. Further, suppose that
 $m(\infty)>0$ and $m(0)>0$. If $\Lambda_m < \Lambda_N \min \bigl(\frac 1{m(\infty)},\frac
1{m(0)}\bigr)$, then there exists a minimizer for $\Lambda_m$.
\end{theorem}

\medskip

\begin{remark}\label{r1}
$\Lambda_m$ has a minimizer also in the following cases, corresponding formally to $\Lambda_+$ or $\Lambda_-$
taking the value $+\infty$.

\medskip

 (i) Let $m(0)=0$ and $m(\infty)>0$. If $\Lambda_m < \frac {\Lambda_N}{m(\infty)}$, then a
minimizer for $\Lambda_1(m)$ exists.

\medskip

(ii)Let $m(0)>0$ and $m(\infty)=0$. If $\Lambda_m < \frac {\Lambda_N}{m(0)}$, then a minimizer for
$\Lambda_1(m)$ exists.

\medskip

(iii) If $m(0)=m(\infty)=0$, $m(x) \geq 0$ and $\not\equiv 0$ on $\RN$, then $\Lambda_m$ has a minimizer.
\end{remark}

\medskip

We point out that Theorem~\ref{t4} and the results described in Remark~\ref{r1} can be deduced from Theorem
1.2 in \cite{TE}. Unlike in paper \cite{TE}, to obtain Theorem~\ref{t4} we avoided the  use of the
concentration - compactness principle.

We now give  examples of weight functions $m$ satisfying conditions of Theorems~\ref{t3} and 2.3. In general,
functions satisfying this condition have large local maxima.

\medskip

\begin{example}\label{e1}
Let
\[
m_A(x)=\left\{
\begin{array}{ll}
m_1(x)\,& \,\,\, \mbox{ for } \, 0<|x|<1,\\
Am_2(x)\,&\,\,\, \mbox{ for }\, 1\leq |x| \leq 2,\\
m_3(x)\, &\,\,\, \mbox{ for }\, 2<|x|,
\end{array}
\right.
\]
where  $A>0$ is a constant to be chosen later and $m_1: \overline {B(0,1)-\{0\}}\to [0,\infty)$, $m_2: (1 \leq
|x| \leq 2) \to [0,\infty)$ and $m_3: \RN\setminus B(0,2)\to [0,\infty)$ are continuous bounded functions
satisfying the following conditions:
 $m_1(x)=0$ for $|x|=1$, $m_2(x)=0$ for $|x|=1$, $m_2(x)=0$ for $|x|=2$, $m_2(x)>0$ for $1<|x|<2$,
$m_3(x)=0$ for $|x|=2$. Further we assume that
\[
m_3(x)=\frac {a+|x_1||x_2|+\ldots +|x_{N-1}||x_N|}{b+|x|^2}
\]
for $|x| \geq R>2$, where $a>0$, $b>0$ and $R$ constants. A function $m_1(x)$ for small $\delta>0$ is given by
\[
m_1(x)=\frac {|x_1|+\ldots +|x_N|}{|x|}
\]
for $0<|x|\leq \delta <1$. We have
\[
\lim_{j \to \infty} m_A(\gamma_+^jx)=\lim_{j \to \infty} \frac {\gamma_+^{-2j}a+|x_1||x_2|+\ldots +
|x_{N-1}||x_N|}{\gamma_+^{-2j}b+|x|^2}= \frac {|x_1||x_2|+\ldots +|x_{N-1}||x_N|}{|x|^2}=m_+(x)
\]
and
\[
\lim_{j \to \infty} m_A(\gamma_-^{-j}x)=\frac {|x_1|+\ldots +|x_N|}{|x|}=m_-(x).
\]

Both limits are uniform.   Since  $m_-$ and $m_+$ are bounded, $\Lambda_-$ and $\Lambda_+$ are positive and
finite. We have
\[
\Lambda_m=\inf_{D^{1,2}(\RN)-\{0\}} \frac {\int_{\RN} |\nabla u|^2\,dx}{\int_{\RN} \frac
{m_A(x)}{|x|^2}u^2\,dx}\leq\frac 1A\inf_{D^{1,2}(\RN)-\{0\}}\frac {\int_{\RN} |\nabla u|^2\,dx}{\int_{1 \leq
|x|\leq 2} \frac {m_2(x)}{|x|^2}u^2\,dx}< \min(\Lambda_-,\Lambda_+)
\]
for $A$ large. By Theorem~\ref{t3}, $\Lambda_m$ with $m=m_A$ has a minimizer.
\end{example}

\medskip

\begin{example}\label{e2}
Consider a sequence  of functions of the form $m_k(x)=BM_k(x)+Af(x), \,\, k=1,2, \ldots$, where $A>0$, $B>0$
are constants and $M_k$ and $f$ are continuous functions satisfying the following conditions:
\begin{description}
\item[(a)] $M_k(0)=1$, $M_k(x)>0$ on $\RN$, $M_k(\infty)=0$ for $k=1,2, \ldots$,
\item[(b)] $M_k(x)=k$ on $1 <|x|<2$ for $k=1,2, \ldots$,
\item[(c)] $f(x) \geq 0$ on $\RN$, $f(0)=0$ and $f(\infty)=1$.
\end{description}
Then $m_k(0)=B$ and $m_k(\infty)=A$ for $k=1,2, \ldots$. We show that for $k$ sufficiently large $m_k$
satisfies the conditions of Theorem~\ref{t4}. Let $u(x)=\exp (-|x|)$ (one can take any other function from
$D^{1,2}(\RN)$ which is $\not\equiv 0$ on $(1<|x|<2)$). Thus
\[
\Lambda_{m_k} \leq \frac {\int_{\RN} |\nabla (\exp(-|x|))|^2\,dx}{\int_{\RN} \frac
{BM_k(x)+Af(x)}{|x|^2}\exp(-2|x|)\,dx} \leq \frac {\int_{\RN} |\nabla (\exp(-|x|))|^2\,dx}{B \int_{\RN} \frac
{M_k(x)}{|x|^2}\exp(-2|x|)\,dx} \to 0,
\]
as $k \to \infty$. So  we can  find $k_{\circ} \geq 1$ so that
\[
\Lambda_{m_k} < \Lambda_N \min \biggl(\frac 1A,\frac 1B\biggr)\,\, \mbox{ for } \,\, k\geq k_{\circ}.
\]
\end{example}

In Proposition~\ref{p1}, below, we described a class of weight functions $m$ satisfying conditions of
Theorem~\ref{t4}.

\begin{proposition}\label{p1}
Let $m \in C(\RN)$. Suppose that $m(x)\geq 0$, $m(0)>0$ and $m(\infty)>0$. Assume that there exists a ball
$B(x_M,r)$ such that $m(x) \geq m(x_M)>0$ for $x \in B(x_M,r)$ and $0 \not\in \overline{B(x_M,r)}$. If
\begin{equation}\label{a18}
\frac {m(0)}{m(x_M)}, \, \frac {m(\infty)}{m(x_M)} < \frac {r^2(N-2)^2}{2(r+|x_M|)^2(N+1)(N+2)}.
\end{equation}
Then $\Lambda_m <\Lambda_N\min\bigl(\frac 1{m(0)},\frac 1{m(\infty)}\bigr)$. (Hence, there exists a minimizer
for $\Lambda_m$.)
\end{proposition}
\pf Let $u \in H_{\circ}^1(B(x_M,r))-\{0\}$. Then
\[
\int_{B(x_M,r)} \frac {m(x)}{|x|^2}u^2\,dx \geq m(x_M) \int_{B(x_M)} \frac {u^2}{|x|^2}\,dx\geq \frac
{m(x_M)}{(r+|x_M|)^2} \int_{B(x_M,r)} u^2\,dx.
\]
Hence
\[
\frac {\int_{B(x_M,r)}|\nabla u|^2\,dx}{\int_{B(x_M,r)} \frac {m(x)}{|x|^2}\,dx} \leq \frac
{(r+|x_M|)^2\int_{B(x_M,r)} |\nabla u|^2\,dx}{m(x_M)\int_{B(x_M,r)} u^2\,dx}.
\]
Since $H_{\circ}^1(B(x_M,r))-\{0\} \subset \{u \in D^{1,2}(\RN);\, \int_{\RN} \frac {m(x)}{|x|^2}u^2\,dx>0\}$
we deduce from the above inequality that
\begin{equation}\label{a19}
\Lambda_m \leq  \frac {(r+|x_M|)^2}{m(x_M)}\lambda_1^D(B(x_M,r)),
\end{equation}
where $\lambda_1^D(B(x_M,r))$ denotes the first eigenvalue for $"-\Delta"$ in $B(x_M,r)$ with the Dirichlet
boundary conditions. We now estimate $\lambda_1^D=\lambda_1^D(B(x_M,r))$. We test $\lambda_1^D$ with
$v(x)=r-|x-x_M|$ for $x \in B(x_M,r)$. We have
\[
\int_{B(x_M,r)} v^2\,dx=\int_{B(0,r)} (r-|x|)^2\,dx=\omega_N \int_0^r(r-s)^2s^{N-1}\,ds=\frac
{2\omega_Nr^{N+2}}{N(N+1)(N+2)}
\]
and
\[
\int_{B(x_M,r)} |\nabla v|^2\,dx=\frac {\omega_Nr^N}N.
\]
Hence
\[
\lambda_1^D \leq \frac {\int_{B(x_M,r)} |\nabla v|^2\,dx}{\int_{B(x_M,r)} v^2\,dx}=\frac {(N+1)(N+2)}{2r^2}.
\]
Combining this with (\ref{a19}) we derive
\[
\Lambda_m \leq \frac {(N+1)(N+2)(r+|x_M|)^2}{2r^2m(x_M)}.
 \]
 Therefore $\Lambda_m <\Lambda_N \min \bigl( \frac 1{m(0)},\frac 1{m(\infty)}\bigr)$ if (\ref{a18}) holds.
 $\hfill\Box$.

 \medskip

 The estimate (\ref{a18}) has terms that are easy to  compute, but are  of course not optimal. In particular, the
 factor $\frac {(N+1)(N+2)}2$can be replaced by the first eigenvalue of the Laplacian on a unit ball with
 Dirichlet boundary conditions.

 If $m(x)$ is a continuous bounded and nonnegative function such that $m(x) \leq m(0)$  on $\RN$ and  $m(0)>0$
(or $m(x)\leq m(\infty)$ on $\RN$, $m(\infty)>0$),   then $\Lambda_m$ does not have a minimizer. Indeed,
suppose that $m(x)\leq m(0)$ on $\RN$ and that $\Lambda_m$ has a minimizer $u$. Then by the Hardy inequality
we obtain
\[
\frac {\Lambda_N}{m(0)}\geq \frac {\int_{\RN} |\nabla u|^2\,dx}{\int_{\RN} \frac {m(x)}{|x|^2}u^2\,dx}\geq
\frac {\int_{\RN}  |\nabla u|^2\,dx}{m(0)\int_{\RN} \frac {u^2}{|x|^2}\,dx} \geq \frac {\Lambda_N}{m(0)}.
\]
So $u$ is a minimizer for $\Lambda_N$, which is impossible.

\medskip

 We now construct a ground state with infinite $D^{1,2}$ norm.

 \medskip

\begin{theorem}\label{t5}
Let $\gamma>1$ and assume that the function $m \in L^{\infty}(\RN)$ satisfies
\begin{equation}\label{a20}
m(\gamma x)=m(x) \,\, \mbox{ for }\,\, x \in \RN.
\end{equation}
Then the form $Q_V$  with $V(x)=\frac {m(x)}{|x|^2}$ and  $\Lambda_V=\Lambda_{\circ}$ (see (\ref{a22}) below)
admits a ground state $v$ satisfying
\begin{equation}\label{a21}
v(\gamma x)=\gamma^{\frac {2-N}2}v(x) \,\, \mbox{ for }\,\, x \in \RN.
\end{equation}
The function  $v$ is uniquely defined by its values on $A_{\gamma}=\{x \in \RN;\, 1<|x|<\gamma\}$ and moreover
the function $v_{|A_{\gamma}}$ is a minimizer for the problem
\begin{equation}\label{a22}
\Lambda_{\circ}=\inf \left\{ \frac {\int_{A_{\gamma}} |\nabla v|^2\,dx}{\int_{A_{\gamma}} \frac
{m(x)}{|x|^2}u^2\,dx}; \, u \in H^1(A_{\gamma})-\{0\},\, u(\gamma x)=\gamma^{\frac {2-N}2}u(x)\,\, \mbox{ \rm
for}\,\, |x|=1\right\}.
\end{equation}
\end{theorem}
\pf The problem (\ref{a22}) is a compact variational problem that has a minimizer $v$ which satisfies the
equation
\[
-\Delta v=\Lambda_{\circ}\frac {m(x)}{|x|^2}v, \,\,\, x \in A_\gamma,
\]
with the Neumann boundary conditions. Since the test functions satisfy $u(\gamma x)=\gamma^{\frac {2-N}2}u(x)$
for $|x|=1$, one has
\begin{equation}\label{a23}
\frac {\partial v}{\partial r}(\gamma x)=\gamma^{-\frac N2}\frac {\partial v}{\partial r}(x) \,\, \mbox{ for
}\,\, |x|=1.
\end{equation}
Note that $|v|$ is also a minimizer, so we may assume that $v$ is nonnegative. We now extend the function $v$
from $A_\gamma$ to $\RN-\{0\}$ by using (\ref{a21}) and denote the extended function again by $v$. Since $v$
satisfies (\ref{a22}), the extended function $v$ is of class $C^1(\RN-\{0\})$ and satisfies the equation
\[
-\Delta v=\Lambda_{\circ}\frac {m(x)}{|x|^2}v
\]
in a weak sense. From this and the Harnack inequality on bounded subsets of $\RN-\{0\}$ it follows that $v$ is
positive on $\RN-\{0\}$ and subsequently there exists a constant $C>0$ such that
\begin{equation}\label{a24}
C^{-1}|x|^{\frac {2-N}2} \leq v(x) \leq C|x|^{\frac {2-N}2}.
\end{equation}
We can now explain  the choice of the exponent $\frac {2-N}2$ in the constraint $u(\gamma x)=\gamma^{\frac
{2-N}2}u(x)$ from  (\ref{a22}): with any other choice the resulting Neumann condition would not yield the
continuity of the derivatives of the extended function $v$ on the spheres $|x|=\gamma^j, \, j \in \mathbb N$. Finally, we show
that $v$ is a ground state for the corresponding quadratic form $Q$ with $V(x)=\Lambda_{\circ}\frac
{m(x)}{|x|^2}$. Using the ground state formula (2.7) from \cite{PTT} and (\ref{a24}), we have  with
$w_k(x)=|x|^{\frac 1k}$ for $|x| \leq 1$ and $w_k(x)=|x|^{-\frac 1k}$ for $|x| \geq 1$,
\begin{eqnarray*}
Q(vw_k)&=&\int_{\RN} v^2|\nabla w_k|^2\,dx \leq C\int_{\RN} |x|^{2-N}|\nabla w_k|^2\,dx\\
&\leq& \frac C{k^2} \int_0^1 r^{-1+\frac 2k}\,dr+\frac C{k^2} \int_1^{\infty} r^{-1-\frac 2k}\,dr\leq \frac Ck
\to 0
\end{eqnarray*}
as $k \to \infty$. Since $vw_k \to v$ uniformly on compact sets, this implies that $v$ is a ground state for
$Q$. By (\ref{a24}) and the Sobolev inequality, $v \not\in D^{1,2}(\RN)$. $\hfill\Box$

\medskip

\section{Perturbations from virtual ground states}
\renewcommand{\theequation}{3.\arabic{equation}}
\setcounter{equation}{0}

\medskip

In this section we show that if a potential term admits a (generalized or \em{large} or \em{virtual}) ground state, then its
weakly continuous perturbations in the suitable direction will admit a ground state with the finite $D^{1,2}$
norm. Then we investigate potentials that do not give rise to a ground state with finite $D^{1,2}$ norm.

\medskip

We need the following existence result.

\medskip

\begin{proposition}\label{p2}
Let $V_{\circ} \in L^{\frac N2,\infty}(\RN)$ be positive on a set of positive measure and let
\begin{equation}\label{a25}
\Lambda_{\circ}=\inf_{u \in D^{1,2}(\RN), \, \int_{\RN} V_{\circ}u^2\,dx=1} \int_{\RN} |\nabla u|^2\,dx.
\end{equation}
Assume that $V_1 \in L^{\frac N2,\infty}(\RN)$ is positive on a set of positive measure and that the
functional $\int_{\RN} \bigl(V_1(x)-V_{\circ}(x)\bigr)u^2\,dx$ is weakly continuous in $D^{1,2}(\RN)$ and let
\begin{equation}\label{a26}
\Lambda_1=\inf_{u \in D^{1,2}(\RN), \, \int_{\RN} V_1u^2\,dx=1} \int_{\RN} |\nabla u|^2\,dx.
\end{equation}
If $\Lambda_1<\Lambda_{\circ}$, then there exists a minimizer for $\Lambda_1$.
\end{proposition}
\pf Let $\{u_k\} \subset D^{1,2}(\RN)$ be a minimizing sequence for (\ref{a26}), that is, $\int_{\RN}
V_1(x)u_k^2\,dx=1$ and $\int_{\RN} |\nabla u_k|^2\,dx \to \Lambda_1$. We may assume that, up to a subsequence,
$u_k \rh w$ in $D^{1,2}(\RN)$ and $L^2(\RN,V_1(x)\,dx)$. Let $v_k=u_k-w$. Then
\begin{eqnarray*}
1&=&\int_{\RN} V_1(x)u_k^2\,dx=\int_{\RN} V_1(x)v_k^2\,dx+\int_{\RN} V_1(x)w^2\,dx+o(1)=\int_{\RN} V_1(x)w^2\,dx\\
&+&\int_{\RN} (V_1(x)-V_{\circ}(x))v_k^2\,dx+\int_{\RN} V_{\circ}(x)v_k^2\,dx+o(1)\\
&=&\int_{\RN} V_{\circ}(x)v_k^2\,dx+\int_{\RN} V_1(x)w^2\,dx+o(1).
\end{eqnarray*}
Let $t=\int_{\RN} V_1(x)w^2\,dx$. Then $\int_{\RN} V_{\circ}(x)v_k^2\,dx \to 1-t$. Assuming that $t<1$ we get
\[
\Lambda_1=\int_{\RN} |\nabla v_k|^2\,dx+\int_{\RN} |\nabla w|^2\,dx+o(1)\geq \Lambda_{\circ}
(1-t)+\Lambda_1t+o(1).
\]
From this we deduce that $\Lambda_1 \geq \Lambda_{\circ}$ which is impossible. Hence $\int_{\RN}
V_1(x)w^2\,dx=1$. From this and the lower semi-continuity of the norm with respect to weak  convergence, we
derive that $w$ is a minimizer and $u_k \to w$ in $D^{1,2}(\RN)$.

\medskip

Proposition~\ref{p2} is related to Theorem 1.7 in \cite{TE} which asserts that a potential of the form
$V(x)=\frac 1{|x|^2}+g(x)$, with a subcritical potential $g$ (for the definition of a subcritical potential
see \cite{TE}) has a  principal eigenfunction. This follows from the fact that $g$ is weakly continuous in
$D^{1,2}(\RN)$ (see \cite{SM}) and the potential $g$ admits a principal eigenfunction.

\medskip

\begin{remark}\label{r2}
\begin{description}
\item[(i)] If $V_1>V_{\circ}$, then $\Lambda_1 \leq \Lambda_{\circ}$, but not necessarily $\Lambda_1 <\Lambda_{\circ}$
\item[(ii)] If  in Proposition~\ref{p2} assumption $\Lambda_1 <\Lambda_{\circ}$ is replaced by $\Lambda_{\circ}
<\Lambda_1$, then $\Lambda_{\circ}$ is attained.
\end{description}
\end{remark}

\medskip

 \begin{example}\label{e4}Let $M$ be a continuous function $\RN$ such that $M\geq 0$, $\not\equiv 0$ on $\RN$ and
$M(0)=M(\infty)=0$.  Define $m_{A,B}(x)=BM(x)+A$, where $A>0$ and $B>0$  are constants. Let $V_1(x)=\frac
{m_{A,B}(x)}{|x|^2}$ and $V_{\circ}(x)=\frac {A}{|x|^2}$. The functional $\int_{\RN}
(V_1(x)-V_{\circ}(x))u^2\,dx=\int_{\RN} \frac {BM(x)}{|x|^2}u^2\,dx$ is weakly continuous in $D^{1,2}(\RN)$.
It is easy to show that for every $A>0$ there exists $B_{\circ}>0$ such that $\Lambda_1 <\Lambda_{\circ}$ for
$B>B_{\circ}$. By Proposition~\ref{p2} $\Lambda_1$ has a minimizer for $B>B_{\circ}$.
\end{example}

\medskip

We now give a sufficient condition for the inequality $\Lambda_1 <\Lambda_{\circ}$.

\begin{theorem}\label{t6}
Suppose that $V_1$ and $V_{\circ}$ satisfy assumptions of Proposition~\ref{p2}. Moreover, assume that the
quadratic form $Q_{V_{\circ}}$ has  a positive ground state $v$, possibly with infinite $D^{1,2}$ norm, and
that, if $\{v_k\} \subset C_{\circ}^{\infty}(\RN)$ is a null sequence corresponding to $\Lambda_{\circ}$, then
\[
\limsup_{k \to \infty} \int_{\RN} (V_1(x)-V_{\circ}(x))v_k^2\,dx>0.
\]
 Then $\Lambda_1 <\Lambda_{\circ}$ and
$\Lambda_1$ has a minimizer.
\end{theorem}
\pf It suffices to show that the inequality
\[
\int_{\RN} |\nabla u|^2\,dx-\Lambda_{\circ} \int_{\RN} V_1(x)u^2\,dx \geq 0
\]
fails for some $u \in D^{1,2}(\RN)$. We have
\begin{eqnarray*}
\int_{\RN} |\nabla v_k|^2\,dx&-&\Lambda_{\circ}\int_{\RN}
V_1(x)v_k^2\,dx=Q_{V_{\circ}}(v_k)-\Lambda_{\circ}\int_{\RN} (V_1(x)-V_{\circ}(x))v_k^2\,dx\\
&=&o(1)-\Lambda_{\circ}\int_{\RN} (V_1(x)-V_{\circ}(x))v_k^2\,dx<0
\end{eqnarray*}
for sufficiently  large $k$, which completes the proof of the theorem. $\hfill\Box$

\medskip

Note that the conditions of Theorem~\ref{t6} are satisfied  if, in particular, $V_1 \geq V_{\circ}$ on $\RN$,
with the strict inequality on a set of positive measure. Indeed, the sequence $\{v_k\}$ converges weakly in
$H_{\mbox{loc}}^1(\RN)$ to $v>0$ and the condition $\limsup_{k \to \infty} \int_{\RN}
(V_1(x)-V_{\circ}(x))v_k^2\,dx>0$ follows from the Fatou lemma.

\medskip

The situation becomes different if $Q_{V_{\circ}}$ does not have a ground state. The absence of the ground
state is stable property under small (in some sense) compact perturbation, but not under compact perturbations
that are not small.

\medskip

\begin{theorem}\label{t7}
Assume that $V_{\circ}$ satisfies the conditions of Proposition~\ref{p2} and that (\ref{a5}) holds (this
occurs under conditions of Theorem~\ref{t2} if $Q_{V{\circ}}$ has no ground state). Let $W$ be as in
(\ref{a5}). Then for every $t \in \bigl(0,\frac 1{\Lambda_{\circ}}\bigr)$ the functional $Q_{V_{\circ}+tW}$
has no ground state and $\Lambda_{V_{\circ}+tW}=\Lambda_{V_{\circ}}$. Furthermore, if the functional
$\int_{\RN} W(x)u^2\,dx$ is weakly continuous in $D^{1,2}(\RN)$, the the same conclusion holds for
$-\infty<t<0$.
\end{theorem}
\pf First we observe that the constants $\Lambda_{\circ}$ and $\Lambda_1$ corresponding to $V_{\circ}$ and
$V_1=V_{\circ}+tW$, respectively, are equal. Indeed, since $V_1>V_{\circ}$, one has $\Lambda_1 \leq
\Lambda_{\circ}$ by monotonicity. On the other hand, it follows from (\ref{a5}) that
\[
\int_{\RN}|\nabla u|^2\,dx-\Lambda_{\circ}\int_{\RN} (V_{\circ}(x)+tW(x))u^2\,dx\geq 0
\]
for $t \in \bigl(0,\frac 1{\Lambda_{\circ}}\bigr)$ which implies $\Lambda_1\geq \Lambda_{\circ}$. Let $v_k \in
C_{\circ}^{\infty}(\RN-Z)$ satisfy $Q_{V_1}(v_k) \to 0$. Then
\[
(1-\Lambda_{\circ}t)\int_{\RN} Wv_k^2\,dx \leq Q_{V_1}(v_k) \to 0,
\]
which implies that,up to subsequence, $v_k \to 0$ a.e. If $v_k$ were a null sequence, it would converge in
$H_{\mbox{loc}}^1(\RN)$ and it would have  a limit zero. Therefore $Q_{V_1}$ admits no null sequence and
consequently no ground state. Assume now that  the functional $\int_{\RN} W(x)u^2\,dx$ is weakly continuous in
$D^{1,2}(\RN)$. Let $\{w_k\}\subset D^{1,2}(\RN)$ be a minimizing sequence for $\Lambda_{\circ}$. If $\{w_k\}$
has  a subsequence weakly convergent in $D^{1,2}(\RN)$ to some $w \neq 0$, then it is easy to see that $|w|$
would be a minimizer for $\Lambda_{\circ}$ and thus a ground state for $Q_{\Lambda_{\circ}}$. Therefore $w_k
\rh 0$. By the weak  continuity of $\int_{\RN} W(x)u^2\,dx$ we get
\[
\int_{\RN} V_1(x)w_k^2\,dx=\int_{\RN} V_{\circ}(x)w_k^2\,dx+o(1)=1+o(1)
\]
and thus
\[
\Lambda_1 \leq \int_{\RN}|\nabla w_k|^2\,dx=\Lambda_{\circ}+o(1).
\]
This yields $\Lambda_1\leq \Lambda_{\circ}$. Then
\begin{eqnarray*}
\int_{\RN} |\nabla u|^2\,dx&-&\Lambda_1 \int_{\RN} V_1(x)u^2\,dx\\
&\geq&\frac {\Lambda_1}{\Lambda_{\circ}}\biggl( \int_{\RN} |\nabla u|^2\,dx-\Lambda_{\circ} \int_{\RN}
V_1(x)u^2\,dx\\
&=&\frac {\Lambda_1}{\Lambda_{\circ}}\biggl(Q_{V_{\circ}}(u)-t\Lambda_{\circ} \int_{\RN} W(x)u^2\,dx\biggr)
\geq \Lambda_1 \int_{\RN} \bigl(\Lambda_{\circ}^{-1}-t\bigr)W(x)u^2\,dx.
\end{eqnarray*}
Since $t<0$, this implies that $Q_{V_1}$ has no ground state. $\hfill\Box$

\medskip

Theorem~\ref{t7} concerns with  small perturbations of a potential that does not change the constant $\Lambda$
or the absence of a ground state. The next theorem shows that a compact perturbation of the potential term
yields a ground state of finite $D^{1,2}(\RN)$ norm.

\medskip

\begin{theorem}\label{t8}
Assume that $V_{\circ}$ satisfies conditions of Proposition~\ref{p2} and that $W \in L^{2,\infty}(\RN)$ is
such that the functional $\int_{\RN} W(x)u^2\,dx$ is weakly continuous in $D^{1,2}(\RN)$. Then for every
$\lambda \in \bigl(0,\Lambda_{\circ}\bigr)$ there exists $\sigma \in \mathbb R$ such that $Q_{V_{\circ}+\sigma
W}$ has  a ground state of finite $D^{1,2}(\RN)$ norm corresponding to the energy constant (\ref{a26}).
\end{theorem}
\pf Assume without loss of generality that $W$ is positive on a set of positive measure. Let
$0<\lambda<\Lambda_{\circ}$ and consider
\[
\sigma=\inf_{u \in D^{1,2}(\RN), \,\int_{\RN} W(x)u^2\,dx=1} \lambda^{-1} \biggl(\int_{\RN} |\nabla
u|^2\,dx-\lambda \int_{\RN} V_{\circ}(x)u^2\,dx\biggr).
\]
Since $\biggl(\int_{\RN} |\nabla u|^2\,dx-\lambda \int_{\RN} V_{\circ}(x)u^2\,dx\biggr)^{\frac 12}$ defines an
equivalent norm on $D^{1,2}(\RN)$ it is easy to show that there exists a minimizer for $\sigma$. It is clear
that this minimizer is also a ground state of $Q_{V_{\circ}+\sigma W}$ corresponding to the optimal constant
$\lambda$. $\hfill\Box$

\medskip

If we assume additionally that $W$ is positive on a set of positive measure, then it is easy to show that
$\sigma$ is a continuous decreasing function of $\lambda$ with $\lim_{\lambda \to 0}\sigma(\lambda)=+\infty$
and  $\sigma_{\circ}=\lim_{\lambda \to \Lambda_{\circ}}\sigma(\lambda)\geq 0$. In particular, if (\ref{a5})
holds with a weight $W_{\circ}$ satisfying $W_{\circ} \geq \alpha  W$, then $\sigma_{\circ} \geq \alpha$. In
other words, given $V_{\circ}$ and $W$ as in Theorem~\ref{t8}, the potential $V_{\circ}+\sigma  W$ admits a
ground state whenever $\sigma\geq \sigma_{\circ}$.

\medskip

For further results of that nature we refer to paper \cite{TE}.

\medskip

\section{Behaviour of a ground state  around $0$ }
\renewcommand{\theequation}{4.\arabic{equation}}
\setcounter{equation}{0}

\medskip

In what follows  we consider the potential of the Hardy type $V(x)=\frac {m(x)}{|x|^2}$, where $m(x)$ is
continuous and $m(0)>0$ and $m(\infty)>0$. The corresponding ground state, if it exists, is denoted by
$\phi_1$, which is chosen to be positive on $\RN$. Obviously the ground state $\phi$ satisfies the equation
\begin{equation}\label{aa1}
\Delta u=\Lambda_m \frac {m(x)}{|x|^2}u \,\, \mbox{ in }\, \RN
\end{equation}
in a weak sense.

\medskip

We need the following extension of the Hardy inequality: let $\Omega \subset \RN$ be a bounded domain and $0
\in \bar \Omega$, then for every $\delta>0$ there exists a constant $A(\delta,\Omega)>0$ such that
\begin{equation}\label{hh}
\int_{\Omega} \frac {u^2}{|x|^2}\,dx \leq \bigl(\frac 1{\Lambda_N}+\delta\bigr) \int_{\Omega} |\nabla
u|^2\,dx+A(\delta,\Omega) \int_{\Omega} u^2\,dx
\end{equation}
for every $u \in H^1(\Omega)$ (see \cite{CHA1}).

\medskip

\begin{proposition}\label{p3}
Let
\[
\Lambda_m < \Lambda_N \min \biggl(\frac 1{m(0)},\frac 1{m(\infty)}\biggr).
\]
Then $\phi_1 \in L^{2^*(1+\delta)}\bigl(B(0,r)\bigr)$ for some $\delta>0$ and $r>0$.
\end{proposition}
\pf Let $\Phi \in C^1(\RN)$ be such that $\Phi(x)=1$ on $B(0,r)$, $\Phi(x)=0$ on $\RN-B(0,2r)$, $0\leq
\Phi(x)\leq 1$ on $\RN$ and $|\nabla \Phi(x)| \leq \frac 2r$. For simplicity we set $\lambda=\Lambda_m$,
$u=\phi_1$. We define $v=\Phi^2 u \min(u,L)^{p-2}=\Phi^2uu_L^{p-2}$, where $L>0$ and $p>2$. Testing
(\ref{aa1}) with $v$, we get
\begin{eqnarray*}
\int_{\RN}\bigl( |\nabla u|^2u_L^{p-2}\Phi^2&+&(p-2)\nabla u \nabla u_Lu_L^{p-2}\Phi^2+2\nabla u \nabla \Phi
uu_L^{p-2}\Phi\bigr)\,dx\\
&=&\lambda \int_{\RN} \frac {m(x)}{|x|^2}u^2u_L^{p-2}\Phi^2\,dx.
\end{eqnarray*}
Applying the Young inequality to the third term on the left side, we get
\begin{eqnarray*}
(1-\eta)\int_{\RN}|\nabla u|^2u_L^{p-2}\Phi^2\,\,dx&+&(p-2) \int_{\RN} \nabla u \nabla
u_Lu_L^{p-2}\Phi^2\,dx\\
&\leq&\lambda\int_{\RN} \frac {m(x)}{|x|^2}u^2u_L^{p-2}\Phi^2\,dx+C(\eta) \int_{\RN} u^2u_L^{p-2}|\nabla
\Phi|^2\,dx,
\end{eqnarray*}
where $\eta>0$ is a small number to be suitably chosen. Since the second integral on the left side is
nonnegative, this inequality can be rewritten in the following form
\begin{eqnarray*}
(1-\eta)\int_{\RN}|\nabla u|^2u_L^{p-2}\Phi^2\,dx&+&(1-\eta)(p-2)\int_{\RN} \nabla u\nabla u_L
u_L^{p-2}\Phi^2\,dx\\
&\leq&\lambda\int_{\RN} \frac {m(x)}{|x|^2}u^2u_L^{p-2}\Phi^2\,dx+C(\eta)\int_{\RN} u^2u_L^{p-2}|\nabla
\Phi|^2\,dx.
\end{eqnarray*}
Multiplying this inequality by $\frac {p+2}4$ and noting that $\frac {p+2}4>1$, we get
\begin{eqnarray}\label{a27}
(1-\eta)\left[ \int_{\RN} |\nabla u|^2u_L^{p-2}\Phi^2\,dx\right.&+&\left.\frac {p^2-4}4 \int_{\RN} \nabla u
\nabla
u_L u_L^{p-2}\Phi^2\,dx\right]\\
&\leq&\frac {\lambda(p+2)}4 \int_{\RN} \frac {m(x)}{|x|^2}u^2u_L^{p-2}\Phi^2\,dx \nonumber\\
&+&\frac {C(\eta)(p+2)}4\int_{\RN} u^2u_L^{p-2}|\nabla \Phi|^2\,dx. \nonumber
\end{eqnarray}
We now observe that
\[
\int_{\RN} |\nabla \bigl(uu_L^{\frac p2-1}\bigr)|^2\Phi^2\,dx=\int_{\RN} |\nabla u|^2u_L^{p-2}\Phi^2\,dx+\frac
{p^2-4}4 \int_{\RN} |\nabla u_L|^2u_L^{p-2}\Phi^2\,dx.
\]
Hence (\ref{a27}) takes the form
\begin{eqnarray}\label{a28}
(1-\eta)\int_{\RN} |\nabla \bigl(uu_L^{\frac p2-1}\bigr)|^2\Phi^2\,dx&\leq&\frac {\lambda(p+2)}4\int_{\RN}
\frac {m(x)}{|x|^2} u^2u_L^{p-2}\Phi^2\,dx\\
&+&\frac {C(\eta)(p+2)}4\int_{\RN} u^2u_L^{p-2}|\nabla \Phi|^2\,dx. \nonumber
\end{eqnarray}
Since $\frac {\lambda m(0)}{\Lambda_N} <1$, we can choose $\epsilon_1>0$ so that $\frac
{\lambda}{\Lambda_N}(m(0)+\epsilon_1)<1$. By the continuity of $m$ there exists $0<r_1<r$ such that $m(x)\leq
m(0)+\epsilon_1$ for $x \in B(0,r_1)$. This is now used to estimate the first integral on the right side of
(\ref{a28}):
\begin{eqnarray*}
\frac {\lambda(p+2)}4 \int_{\RN} \frac {m(x)}{|x|^2}u^2u_L^{p-2}\Phi^2\,dx& \leq& \frac {\lambda(p+2)}4
\int_{B(0,r_1)} \frac {m(0)+\epsilon_1}{|x|^2}u^2u_L^{p-2}\,dx\\
&+&\frac {\lambda(p+2)\|m\|_{\infty}}{4r_1^2} \int_{B(0,2r)} u^2u_L^{p-2}\,dx.
\end{eqnarray*}
Applying the Hardy inequality (\ref{hh}), we get
\begin{eqnarray*}
\frac {\lambda(p+2)}4 \int_{B(0,r)}\frac {m(x)}{|x|^2}u^2u_L^{p-2}\,dx&\leq&\frac
{\lambda(p+2)}4(m(0)+\epsilon_1)\bigl(\frac 1{\Lambda_N}+\epsilon\bigr)\int_{B(0,r_1)} |\nabla
\bigl(uu_L^{\frac p2-1})|^2\,dx\\
&+&\biggl(\frac {\lambda (p+2)}4 A(B(0,r_1),\epsilon)\\
&+&\frac {\lambda(p+2)\|m\|_{\infty}}{4r_1^2}\biggr)\int_{B(0,2r)} \bigl(uu_L^{\frac p2-1}\bigr)^2\,dx
\end{eqnarray*}
for every $\epsilon>0$. Inserting this estimate into (\ref{a28}) we obtain
\begin{equation}\label{a29}
\biggl(1-\eta-\frac {\lambda(p+2)}4(m(0)+\epsilon_1)\bigl(\frac 1{\Lambda_N}+\epsilon\bigl)\biggr)
\int_{B(0,r)} |\nabla \bigl(uu_L^{\frac p2-1}\bigr)|^2\,dx\leq C_1\int_{B(0,2r)} \bigl(uu_L^{\frac
p2-1}\bigr)^2\,dx,
\end{equation}
where $C_1=\frac {\lambda (p+2)}4 A(B(0,r_1),\epsilon)+\frac {\lambda(p+2)\|m\|_{\infty}}{4r_1^2}+\frac
{(p+2)C(\eta)}{r^2}$. We put $p=2+\delta$, $\delta>0$. We now observe that we can choose $\delta$ and
$\epsilon$ so small that
\[
\lambda\bigl(1+\frac {\delta}4\bigr)(m(0)+\epsilon_1)\bigl(\frac 1{\Lambda_N}+\epsilon\bigr)=\frac
{\lambda}{\Lambda_N}\bigl(1+\frac {\delta}4\bigr)(m(0)+\epsilon_1)+\lambda \epsilon\bigl(1+\frac
{\delta}4\bigr)(m(0)+\epsilon_1)<1.
\]
We point out that we have used here the inequality $\frac {\lambda}{\Lambda_N}(m(0)+\epsilon_1)<1$. With this
choice of $\epsilon$ and $\delta$ we now  choose $\eta>0$ so small that
\[
C_2:=1-\eta-\lambda\bigl(1+\frac {\delta}4\bigr)(m(0)+\epsilon_1)\bigl(\frac 1{\Lambda_N}+\epsilon\bigr)>0.
\]
Finally, we apply the Sobolev inequality in $H^1(B(0,r))$ and  deduce
\[
SC_2\biggl(\int_{B(0,r)} |uu_L^{\frac p2-1}|^{2^*}\,dx\biggr)^{\frac 2{2^*}}\leq (C_1+C_2)\int_{B(0,2r)}
\bigl(uu_L^{\frac p2-1}\bigr)^2\,dx,
\]
where $S$ denotes the best Sobolev constant of the embedding of $H^1(B(0,r))$ into $L^{2^*}(B(0,r))$. Letting
$L \to \infty$ we deduce that $u \in L^{2^*(1+\frac {\delta}2)}(B(0,r))$. So the assertion holds with
$\delta_{\circ}=\frac {\delta}2$.  $\hfill\Box$

\medskip

We now establish the higher integrability property of the principal eigenfunction on $\RN\setminus B(0,R)$.
Although this  will not be used in the sequel, we add it for the sake of completness. We denote by
$D^{1,2}(\RN\setminus B(0,R))$ the Sobolev space defined by
\[
D^{1,2}(\RN\setminus B(0,R))=\{u: \, \nabla u \in L^2(\RN\setminus B(0,R)) \, \mbox{ and }\, u \in
L^{2^*}(\RN\setminus B(0,R))\}.
\]

\medskip

\begin{lemma}\label{l2}
For every $\delta>0$, there exists a constant  $A=A(\delta,R)>0$ such that
\[
\int_{|x| \geq R} \frac {u^2}{|x|^2}\,dx \leq \bigl(\frac 1{\Lambda_N}+\delta\bigr) \int_{|x|\geq R} |\nabla
u|^2\,dx+A\int_{R \leq |x|\leq R+1} u^2\,dx
\]
for every $u \in D^{1,2}(\RN\setminus B(0,R))$.
\end{lemma}
\pf Let $\Phi \in C^1(\RN)$ be such that $\Phi(x)=0$ on $\overline{B(0,R)}$, $\Phi(x)=1$ on $\RN\setminus
B(0,R+1)$, $0\leq \Phi(x)\leq 1$ on $\RN\setminus B(0,R)$ and $|\nabla \Phi(x)| \leq \frac  2R$ on $\RN$. Then
$u\Phi \in D^{1,2}(\RN)$ and by the Hardy and Young inequalities, we have
\begin{eqnarray*}
\int_{|x| \geq R} \frac {u^2}{|x|^2}\,dx&=&\int_{|x|\geq R} \frac {(u\Phi)^2}{|x|^2}\,dx+\int_{|x|\geq R}
\frac {(1-\Phi^2)u^2}{|x|^2}\,dx\\
&\leq&\Lambda_N^{-1}\int_{|x|\geq R} |\nabla (u\Phi)|^2\,dx+\frac 1{R^2}\int_{R\leq |x| \leq R+1} u^2\,dx\\
&\leq&\Lambda_N^{-1} \int_{|x|\geq R} |\nabla u|^2\Phi^2\,dx +\Lambda_N^{-1} \int_{|x| \geq R} u^2|\nabla
\Phi|^2\,dx\\
&+&2\Lambda_N^{-1} \int_{|x|\geq R} u\Phi\nabla u\nabla \Phi\,dx
+\frac 1{R^2} \int_{R \leq |x| \leq R+1} u^2\,dx\\
&\leq&\bigl(\Lambda_N^{-1}+\delta\bigr) \int_{|x| \geq R} |\nabla u|^2\,dx+\bigl( \Lambda_N^{-1}+C(\delta)
\bigr) \int_{|x| \geq R} u^2|\nabla \Phi|^2\,dx\\
 &+&\frac 1{R^2} \int_{R\leq |x|\leq R+1} u^2\,dx
\end{eqnarray*}
and the result follows with $A(\delta,R)=\frac 4{R^2}\bigl(\Lambda_N^{-1}+C(\delta)\bigr)+\frac 1{R^2}$.
$\hfill\Box$

\medskip

\begin{proposition}\label{p4}
Suppose that  $m(\infty)>0$ and $\Lambda_m<\Lambda_N\min\biggl(\frac 1{m(0)},\frac 1{m(\infty)}\biggr)$. Let
$\phi_1$ be the principal eigenfunction of problem (\ref{aa1}). Then there exist $\delta>0$ and $R>0$ such
that $\phi \in L^{2^*(1+\delta)}(\RN\setminus B(0,R))$.
\end{proposition}
\pf We modify the argument used in the proof of Proposition~\ref{p3}. Since $\Lambda_m <\frac
{\Lambda_N}{m(\infty)}$, there exist $\epsilon>0$ and $R>0$ such that $\frac {\Lambda_m}{\Lambda_N}
(m(\infty)+\epsilon)<1$ and $m(x)<m(\infty)+\epsilon$  for $|x| \geq R$. Let $\Psi \in C^1(\RN)$ be such that
$\Psi(x)=0$ on $B(0,R)$, $\Psi(x)=1$ on $\RN-B(0,R+1)$, $0 \leq \Psi(x) \leq 1$ on $\RN$ and $|\nabla \Psi(x)|
\leq \frac 2{R}$ on $\RN$. Let $\lambda =\Lambda_m$, $u=\phi_1$ and $v=uu_L^{p-2}\Psi^2$, where $L>1$, $p>2$
and $u_L=\min(u,L)$. It is clear that $v \in D^{1,2}(\RN)$. Testing  (\ref{aa1}) with $v$ and applying the
Young inequality, we obtain
\begin{eqnarray*}
(1-\eta) \int_{\RN} |\nabla u|^2u_L^{p-2}\Psi^2\,dx&+&(p-2)\int_{\RN} \nabla u \nabla u_Lu_L^{P-2}\Psi^2\,dx\\
&\leq&\lambda \int_{\RN} \frac {m(x)}{|x|^2}u^2u_L^{p-2}\Psi^2\,dx+C(\eta) \int_{\RN} u^2u_L^{p-2}|\nabla
\Psi||^2\,dx.
\end{eqnarray*}
From this, as in the proof of Proposition~\ref{p3}, we derive that
\begin{eqnarray}\label{b1}
(1-\eta) \int_{\RN} |\nabla\bigl(uu_L^{\frac p2-1}\bigr)|^2\Psi^2\,dx&\leq&\frac {\lambda (p+2)}4 \int_{\RN}
\frac
{m(x)}{|x|^2}u^2u_L^{p-2}\Psi^2\,dx\\
&+&\frac {C(\eta)(p+2)}4\int_{\RN} u^2u_L^{p-2}|\nabla \Psi|^2\,dx. \nonumber
\end{eqnarray}
We now estimate the first integral on the right side of (\ref{b1}). Using Lemma~\ref{l2} we have for every
$\epsilon_1>0$
\begin{eqnarray*}
\int_{\RN} \frac {m(x)}{|x|^2}u^2u_L^{p-2}\Psi^2\,dx &\leq& (m(\infty)+\epsilon) \int_{|x|\geq R+1} \frac
{\bigl(uu_L^{\frac p2-1}\bigr)^2}{|x|^2}\,dx\\
&+&(m(\infty)+\epsilon) \int_{R \leq |x|\leq R+1} \frac {\bigl(uu_L^{\frac p2-1}\bigr)^2}{|x|^2}\,dx\\
&\leq &\bigl(\Lambda_N^{-1}+\epsilon_1\bigr)(m(\infty)+\epsilon)\int_{|x|\geq R+1}
|\nabla \bigl(uu_L^{\frac p2-1}\bigr)|^2\,dx\\
&+&A(\epsilon_1,R)(m(\infty)+\epsilon) \int_{R+1\leq |x| \leq R+2} \bigl(uu_L^{\frac p2-1}\bigr)^2\,dx\\
&+&\frac {m(\infty) +\epsilon)}{R^2}\int_{R \leq |x|\leq R+1} \bigl(uu_L^{\frac p2-1}\bigr)^2\,dx.
\end{eqnarray*}
Inserting this into (\ref{b1}) we obtain
\begin{eqnarray}\label{b2}
&&\left[1-\eta -\frac {\lambda (p+2)}4 \bigl(\Lambda_N^{-1} +\epsilon_1\bigr)(m(\infty)+\epsilon)\right]
\int_{|x|\geq R+1} |\nabla\bigl(uu_L^{\frac p2-1}\bigr)|^2\,dx\\
&\leq&C_1(\delta,\epsilon_1,R)\int_{R \leq |x|\leq R+2} \bigl(uu_L^{\frac p2-1}\bigr)^2\,dx,\nonumber
\end{eqnarray}
where
\[
C_1(\delta,\epsilon_1,R):=\frac {\lambda (p+2)}4(m(\infty)+\epsilon)A(\epsilon_1,R)+\frac {\lambda
(p+2)}{4R^2} (m(\infty)+\epsilon)+\frac {C(\eta)(p+2)}{R^2}.
\]
We now set $p=2+\delta$. We choose $\delta>0$ and $\epsilon_1>0$ such that
\[
\lambda (1+\frac {\delta}4)\bigl(\Lambda_N^{-1}+\epsilon_1\bigr)(m(\infty)+\epsilon)<1.
\]
Then we choose $\eta>0$ small enough to guarantee the inequality
\[
C_2:=1-\eta-\lambda \bigl(1+\frac {\delta}4\bigr)\bigl(\Lambda_N^{-1}+\epsilon_1\bigr)(m(\infty)+\epsilon)>0.
\]
Having chosen $\epsilon_1$ and $\delta$ we apply the Sobolev inequality to deduce from (\ref{b2})
\[
SC_2\biggl( \int_{|x| \geq R+1} |\bigl(uu_L^{\frac p2-1}\bigr)|^{2^*}\,dx\biggr)^{\frac 2{2^*}} \leq
C_1\int_{R\leq |x| \leq R+1} \bigl(uu_L^{\frac p2-1}\bigr)^2\,dx,
\]
where $S$ is the best Sobolev constant for the embedding of $D^{1,2}(\RN-B(0,R+1))$ into
$L^{2^*}(\RN-B(0,R+1))$. Letting $L \to \infty$, the result follows. $\hfill\Box$

\medskip

Continuing with the above notations $\lambda=\Lambda_m$, $u=\phi_1$, we put $u=|x|^{-s}v$, with $s>0$ to be
chosen later. We have
\[
\mbox{ div }\bigl(|x|^{-2s}\nabla v\bigr)=-\lambda
|x|^{-2-s}m(x)u+u\bigl(-s^2|x|^{-s-2}+sN|x|^{-s-2}-2s|x|^{-s-2}\bigr).
\]
We now consider the above equation in a small ball $B(0,r)$. Since
 \[
\lambda=\Lambda_m<\Lambda_N\min \biggl(\frac 1{m(0)},\frac 1{m(\infty)}\biggr)\leq \frac {\Lambda_N}{m(0)},
\]
there exists $r>0$ (small enough) such that $\lambda \max_{x \in B(0,r)} m(x) <\Lambda_N$. Let
$s=\sqrt{\Lambda_N}-\sqrt{\Lambda_N-\lambda \bar m_r}$ with $\bar m_r=\max_{x \in B(0,r)} m(x)$, then
\begin{equation}\label{a30}
-\mbox{div }\bigl(|x|^{-2s}\nabla v\bigr)\leq 0 \,\,\, \mbox{ in }\,\, B(0,r).
\end{equation}
Let $\underbar m_r=\min_{x \in B(0,r)} m(x) $ and set $s=\sqrt{\Lambda_N}-\sqrt{\Lambda_N-\lambda\underbar
m_r}$. Then
\begin{equation}\label{a31}
-\mbox{div }\bigl(|x|^{-2s}\nabla v\bigr)\geq 0 \,\,\, \mbox{ in }\,\, B(0,r).
\end{equation}

\medskip

\begin{proposition}\label{p5}
Let $m(0)>0$ and
\[
\Lambda_m<\Lambda_N\min \biggl(\frac 1{m(0)},\frac 1{m(\infty)}\biggr).
 \]
 Then there exists $r>0$ such that
\begin{equation}\label{a32}
M_1|x|^{-(\sqrt{\Lambda_N}-\sqrt{\Lambda_N-\lambda\underbar m_r})} \leq \phi_1(x)\leq M_2
|x|^{-(\sqrt{\Lambda_N}-\sqrt{\Lambda_N-\lambda_m \bar m_r})}
\end{equation}
for $x \in B(0,r)$ and some constants $M_1>0$, $M_2>0$.
\end{proposition}
The lower bound follows from Proposition 2.2 in \cite{CHE}. To apply it we need inequality (\ref{a31}). To
establish the upper bound,  we modify the argument used in paper \cite{HA}. Let $\eta$ be a $C^1$ function
such that $\eta(x)=1$ on $B(0,r)$, $\eta(x)=0$ on $\RN\setminus B(0,\rho)$ and $|\nabla \eta(x)|\leq \frac
2{\rho-r}$ on $\RN$,  where $0<r<\rho$. We use as a test function in (\ref{a30})
$w=\eta^2vv_l^{2(t-1)}=\eta^2v\min(v,l)^{2(t-1)}$, where $l,t>1$. Substituting into (\ref{a30}), we obtain
\begin{equation}\label{a33}
\int_{\RN} |x|^{-2s}\bigl(2\eta vv_l^{2(t-1)}\nabla v\nabla \eta+\eta^2v_l^{2(t-1)}|\nabla
v|^2+2(t-1)\eta^2v_l^{2(t-1)}|\nabla v_l|^2\bigr)\,dx \leq 0,
\end{equation}
where $s=\sqrt{\Lambda_N}-\sqrt{\Lambda_N-\lambda \bar m_r}$. By the Young inequality, for every $\epsilon>0$
there exists $C(\epsilon)>0$ such that
\begin{eqnarray*}
2\int_{\RN} |x|^{-2s}\eta vv_l^{2(t-1)}\nabla \eta \nabla v\,dx&\leq&\epsilon \int_{\RN}
|x|^{-2s}\eta^2v_l^{2(t-1)}|\nabla v|^2\,dx\\
&+&C(\epsilon) \int_{\RN} |x|^{-2s}|\nabla \eta|^2v^2v_l^{2(t-1)}\,dx.
\end{eqnarray*}
Taking $\epsilon=\frac 12$, we derive from (\ref{a33}) that
\begin{eqnarray}\label{a34}
\int_{\RN} |x|^{-2s}\bigl(\eta^2v_l^{2(t-1)}|\nabla v|^2&+&2(t-1)\eta^2v_l^{2(t-1)}|\nabla v_l|^2\bigr)\,dx\\
&\leq&C \int_{\RN} |x|^{-2s}|\nabla \eta|^2v^2v_l^{2(t-1)}\,dx,\nonumber
\end{eqnarray}
where $C>0$ is a constant independent of $l$. To proceed further we use the Caffarelli - Kohn - Nirenberg
inequality \cite{CKN}:
\begin{equation}\label{a35}
\biggl( \int_{B(0,\rho)} |x|^{-bp}|w|^p\,dx\biggr)^{\frac 2p}\leq C_{a,b} \int_{B(0,\rho)} |x|^{-2a}|\nabla
w|^2\,dx
\end{equation}
for every $w \in H_{\circ}^1\bigl(B(0,\rho),|x|^{-2a}\,dx\bigr)$, where $-\infty <a<\frac {N-2}2$, $a\leq b
\leq a+1$, $p=\frac {2N}{N-2+2(b-a)}$ and $C_{a,b}>0$ is a constant depending on $a$ and $b$. We choose
\[
a=b=\sqrt{\Lambda_N}-\sqrt{\Lambda_N-\lambda \bar m_r}< \frac {N-2}2.
\]
In this case we have $p=2^*$. We then deduce from (\ref{a34}) and (\ref{a35}) with $w=\eta vv_l^{t-1}$, that
\begin{eqnarray}\label{a36}
\biggl(\int_{\RN} |x|^{-2^*s}|\eta vv_l^{t-1}|^{2^*}\,dx\biggr)^{\frac 2{2^*}}&\leq& C_{a,b} \int_{\RN}
|x|^{-2s}|\nabla \bigl(\eta vv_l^{t-1}\bigr)|^2\,dx\\
&\leq&2C_{a,b} \int_{\RN} |x|^{-2s}\bigl(|\nabla \eta|^2v^2v_l^{2(t-1)}+\eta^2v_l^{2(t-1)}|\nabla
v|^2\nonumber\\
&+&(t-1)^2\eta^2v_l^{2(t-1)}|\nabla v_l|^2\bigr)\,dx\nonumber\\
&\leq&Ct \int_{\RN} |x|^{-2^*s}|\nabla \eta|^2v^2v_l^{2(t-1)}\,dx. \nonumber
\end{eqnarray}
We now observe that
\[
 \int_{\RN} |x|^{-2^*s}|\eta|^{2^*}v^2v_l^{2^*t-2}\,dx \leq \int_{\RN} |x|^{-2^*s}|\eta vv_l^{t-1}|^{2^*}\,dx.
\]
Indeed, to show this we need to check that $v^2v_l^{2^*t-2}\leq v_l^{2^*(t-1)}v^{2^*}$ on $\mbox{ supp }
\eta$. This can be verified by considering the cases $v_l=l$ and $v_l=v$. The above inequality allows us to
rewrite (\ref{a36}) as
\[
\biggl( \int_{\RN} |x|^{-2^*s}|\eta|^{2^*}v^2v_l^{2^*t-2}\,dx\biggr)^{\frac 2{2^*}}\leq Ct\int_{\RN}
|x|^{-2^*s}|\nabla \eta|^2v^2v_l^{2(t-1)}\,dx.
\]
Due to the properties of the function $\eta$, the above inequality becomes
\begin{equation}\label{a37}
\biggl( \int_{B(0,r)} |x|^{-2^*s}v^2v_l^{2^*t-2}\,dx\biggr)^{\frac 2{2^*}}\leq \frac
{Ct}{(\rho-r)^2}\int_{B(0,\rho)} |x|^{-2^*s}v^2v_l^{2(t-1)}\,dx.
\end{equation}
One can easily check that the resulting integral on the right side is of (\ref{a37}) is finite. We now choose
$\frac N{N-2}<t^*<(1+\delta_{\circ})\frac N{N-2}$, where $\delta_{\circ}$ is a constant from
Proposition~\ref{p3}. We define the sequence $t_j=t^*\bigl(\frac {2^*}2\bigr)^j$, $j=0,1, \ldots$. Setting
$t=t_j$ in (\ref{a37}), we obtain
\[
\biggl( \int_{B(0,r)} |x|^{-2^*s}v^2v_l^{2t_{j+1}-2}\,dx\biggr)^{\frac 1{2t_{j+1}}}\leq \biggl(\frac {
Ct_j}{(\rho-r)^2}\biggr)^{\frac 1{2t_j}}\biggl(\int_{B(0,\rho)} |x|^{-2^*s}v^2v_l^{2t_j-2}\,dx\biggr)^{\frac
1{2t_j}}.
\]
We put $r_j=\rho_{\circ}\bigl(1+\rho_{\circ}^j\bigr)$, $j=0,1, \ldots$ with $\rho_{\circ}$ small. Substituting
in the last inequality $\rho=r_j$, $r=r_{j+1}$, we obtain
\begin{equation}\label{a38}
\biggl( \int_{B(0,r_{j+1})} |x|^{-2^*s}v^2v_l^{2t_{j+1}-2}\,dx\biggr)^{\frac 1{2t_{j+1}}}\leq \biggl(\frac
{Ct_j}{(\rho_{\circ}-\rho_{\circ}^2)^2\rho_{\circ}^{2j}}\biggr)^{\frac
1{2t_j}}\biggl(\int_{B(0,r_j)}|x|^{-2^*s}v^2v_l^{2t_j-2}\,dx\biggr)^{\frac 1{2t_j}}.
\end{equation}
Iterating gives
\begin{eqnarray}\label{a39}
\biggl(\int_{B(0,r_{j+1})}&& |x|^{-2^*s}v^2v_l^{2t_{j+1}-2}\,dx\biggr)^{\frac 1{2t_{j+1}}}\\
&\leq& \biggl(\frac C{\rho_{\circ}-\rho_{\circ}^2}\biggr)^{\sum_{j=0}^{\infty}\frac
1{t_j}}\rho_{\circ}^{-\sum_{j=0}^{\infty}\frac 1{t_j}}\prod_{j=0}^{\infty} t_j^{\frac
1{2t_j}}\biggl(\int_{B(0,r_{\circ})} |x|^{-2^*s}v^2v_l^{2t^*-2}\,dx\biggr)^{\frac 1{2^*}}.
\end{eqnarray}
We now notice that infinite sums and the infinite product in the above inequality are finite. Since
$2^*<2t^*<(1+\delta_{\circ})2^*$, we have
\begin{equation}\label{a40}
\int_{B(0,r_{\circ})} |x|^{-2^*s}v^2v_l^{2t^*-2}\,dx \leq \int_{B(0,r_{\circ})}
|x|^{(2t^*-2^*)s}|u|^{2t^*}\,dx \leq r_{\circ}^{(2t^*-2^*)s}\int_{B(0,r_{\circ})}|u|^{2^*t^*}\,dx<\infty.
\end{equation}
We now deduce from (\ref{a39}) and (\ref{a40}) that
\begin{eqnarray*}
\|v_l\|_{L^{2t_{j+1}}(B(0,\rho_{\circ}))}&\leq&\|v_l\|_{L^{2t_{j+1}}(B(0,r_{j+1}))}\\
&\leq&r_{\circ}^{\frac {2^*s}{2t_{j+1}}}\biggl(\int_{B(0,r_{j+1})}
|x|^{-2^*s}v^2v_l^{t_{j+1}-2}\,dx\biggl)^{\frac 1{2t_{j+1}}}\leq C,
\end{eqnarray*}
where $C>0$ is a constant independent of $l$ and $j$. Letting $t_j \to \infty$ we get
$\|v_l\|_{L^{\infty}(B(0,\rho_{\circ}))} \leq C$. Finally, if $l \to \infty$ we obtain
$\|v\|_{L^{\infty}(B(0,\rho_{\circ}))} \leq C$ and this completes the proof. $\hfill\Box$

\end{document}